\documentclass[a4paper]{article}
\usepackage{amsthm,amsmath,amssymb}
\usepackage{enumerate,enumitem}
\usepackage{graphicx}
\usepackage{xfrac}
\usepackage{tabularx,booktabs,multirow,cleveref}
\newcolumntype{C}{>{\centering\arraybackslash}X}
\newcolumntype{D}{>{\centering\arraybackslash}X}

\usepackage{thmtools}
\usepackage{thm-restate}

\DeclareGraphicsExtensions{.pdf,.png,.eps}
\graphicspath{{figs/}}

\newtheorem{theorem}{Theorem}
\newtheorem{lemma}{Lemma}

\newtheorem{corollary}[theorem]{Corollary}

\newtheorem*{claim*}{Claim}

\theoremstyle{remark}

\newcommand{\cC}{\ensuremath{\mathcal{C}}}
\newcommand{\G}{\ensuremath{\mathcal{G}}}

\newcommand{\cH}{\ensuremath{\mathcal{H}}}

\newcommand{\ex}{{\rm  ex}}
\newcommand{\de}{\mathbf{d}}

\newcommand{\abs}[1]{\left\lvert{#1}\right\rvert}
\newcommand{\floor}[1]{\left\lfloor{#1}\right\rfloor}
\newcommand{\ceil}[1]{\left\lceil{#1}\right\rceil}

\newcommand{\ones}{\ensuremath{\textbf{1}}}

\usepackage[usenames,dvipsnames]{xcolor}

\newcommand\ordstar{\mathord{\star}}

\begin{document}

\title{Generalized Tur\'an densities in the hypercube}
\author{
 Maria Axenovich\thanks{Karlsruhe Institute of Technology, Karlsruhe, Germany, \texttt{maria.aksenovich@kit.edu}.}
 \and Laurin Benz\thanks{Karlsruhe Institute of Technology, Karlsruhe, Germany, \texttt{laurin.benz@kit.edu}.}
 \and David Offner\thanks{Carnegie Mellon University, Pittsburgh, PA, USA, \texttt{doffner@andrew.cmu.edu}.}
 \and Casey Tompkins\thanks{Karlsruhe Institute of Technology, Karlsruhe, Germany and Institute for Basic Science, Daejeon, South
  \indent \hspace{1.5mm} Korea and Alfr\'ed R\'enyi  Institute of  Mathematics, Budapest, Hungary
\texttt{ctompkins496@gmail.com}.}
 }
\date{\today}

\maketitle

\begin{abstract}
A classical extremal, or Tur\'an-type problem asks to determine $\ex (G, H)$, the largest number of edges in a subgraph of a graph $G$ which does not contain a subgraph isomorphic to $H$.  Alon and Shikhelman introduced the so-called generalized extremal number $\ex(G,T,H)$, defined to be the maximum number of  subgraphs isomorphic to $T$ in a subgraph of $G$ that contains no subgraphs isomorphic to $H$.  In this paper we investigate the case when $G = Q_n$, the hypercube of dimension $n$, and $T$ and $H$ are smaller hypercubes or cycles.
\end{abstract}

\section{Introduction}

For a graph $G$ we write $||G||$ for the number of edges of $G$ and $E(G)$ for the set of edges of~$G$.
For two graphs $G$ and $G'$, we write that $G\cong G'$ if the graphs are isomorphic, and we write $G'\subseteq G$ if $G'$ is a subgraph of $G$.
Given a graph $T$, we refer to each $G' \subseteq G$ isomorphic to $T$ as a {\it copy} of $T$.
If a graph $G$ does not have a subgraph isomorphic to $H$, we say that $G$ is $H$-{\it free}.
We denote a complete graph on $n$ vertices by $K_n$, a cycle on $n$ vertices by $C_n$, and a hypercube of dimension $n$ by $Q_n$.
Recall that a hypercube of dimension $n$ is a graph whose vertices are binary sequences of length $n$ and whose edges are pairs of vertices that differ in exactly one position, i.e., having Hamming distance one.
For a graph $G$ and its subgraphs $T$ and $H$, let
\[ N(G,T) = \abs{\{ T':  T'\subseteq G, T'\cong T\}}, \]
\[ \ex(G,T,H) = \max \{ N(G',T) \mid  G'\subseteq G,  G' \mbox{ is $H$-free}\}, \]
\[ \de(G,T,H) =\frac{ \ex(G,T,H)}{N(G,T)}. \]

In a plain language, $N(G,T)$ counts the number of copies of $T$ in $G$,
$\ex(G,T,H)$ is the maximum number of copies of $T$ in an $H$-free subgraph of $G$ and  $\de(G,T,H)$ gives the largest proportion of the number of copies of $T$ in an $H$-free subgraph of $G$.
We refer to $\ex(G,T, H)$ as the {\it generalized extremal function for $H$ in $G$ with respect to $T$} and  $\de(G,T,H) $ as the {\it generalized Tur\'an density of $H$ in $G$ with respect to $T$}.
When $T=K_2$, $\ex(G,T,H)$  is equal to the classical extremal function $\ex(G,H)$ counting the maximum number of edges in an $H$-free subgraph of $G$. In particular, when $G=K_n$,  $\ex(G,H)= \ex(K_n,H)= \ex(n,H)$.
Extremal functions of graphs have been studied extensively.  Since we are concerned here with the ground graph $G=Q_n$,  we will provide a summary of known extremal functions $\ex(Q_n,H)$ in \Cref{extremal-hypercube}.
\vspace{\baselineskip}

In the context of generalized extremal functions, mostly the case $G=K_n$ has been considered.  Already in 1949, Zykov~\cite{zykov}  (and later independently Erd\H{o}s~\cite{erdos}) determined the value of $\ex(n,K_r,K_t)$ for all $r$ and $t$, thereby generalizing the classical theorem of Tur\'an~\cite{turan}. Subsequently, other pairs of graphs have been considered. Particular attention was paid to determining $\ex(n,C_5,C_3)$. The value was estimated within a factor of $1.03$ by Gy\H{o}ri \cite{gyori} and later determined exactly through the method of flag algebras by Hatami {et al.}~\cite{hatami} and independently Grzesik \cite{grzesik}. The systematic study of the function $\ex(n,T,H)$ was initiated by Alon and Shikhelman~\cite{alon}. The problem of determining $\ex(G(n,p),T,H)$ for the random graph $G(n,p)$ was recently investigated by Samotij and Shikhelman \cite{samotij} (See~\cite{alon0} for the case when $T$ is a clique).  Furthermore, Alon and Shikhelman \cite{alon2} recently proved some algorithmic properties of the generalized extremal function $\ex(G,T,H)$.
\vspace{\baselineskip}

In this paper we prove the following statements about $\ex(Q_n, T, H)$, where $T$ and $H$ are either cycles or smaller hypercubes. For asymptotic Landau notations $o, O$, etc., we shall always consider $n$ tending to infinity while the other parameters are fixed and our terms $o(f(n))$ are assumed to be non-negative.

\begin{restatable}{theorem}{QnQlQk} \label{thm:QnQlQk}
    For any integers $k$ and $\ell$ with $2 \leq \ell < k$ and sufficiently large integer $n$, we have
    \[  \max \left\{ 1 - \frac{\ell}{k}, 1- \frac{4 \binom{\ell+2}{3}}{k(k+2)}\right\}  \leq  \de(Q_n, Q_\ell, Q_k) \leq
    \min \left\{ 1-\frac{\ell 2^{\ell}}{k2^k},  1- \alpha \frac{\log k}{k2^k}\right\}, \]
    for a positive constant $\alpha$.
\end{restatable}

The second expression in the lower bound is larger  whenever $k > 2 \ell^2 /3 + 2\ell - 2/3$.

\begin{restatable}{theorem}{QnCfourCsix} \label{thm:QnC4C6}
    For any sufficiently large integer $n$,~
    $ 0.25  n^{-1}  \leq \de(Q_n, C_4, C_6) \leq 0.36578 n^{-1}.$
\end{restatable}

Note that an edge in $Q_n$ corresponds to two binary vectors that differ in exactly one position that is referred to as a {\it star} or a flip position.  
Since we do not have a precise expression for $N(Q_n, C_{2\ell})$, we use an asymptotic result for the theorems about counting $C_{2 \ell}$'s. Let $z_{k,\ell}$ be the number of $C_{2 \ell}$'s in a $Q_k$ using exactly  $k$  different  star positions on its edges.  For a formal definition, see Section \ref{extremal-hypercube}.

\begin{restatable}{theorem}{QnClCsix} \label{thm:QnC2lC6}
    For any integer $\ell$ with $\ell \geq 4$ and sufficiently large integer $n$, we have
    \[\left(4^{\ell+1} z_{\ell,\ell} (1+o(1)) \right)^{-1} \leq  \de(Q_n, C_{2 \ell}, C_6) \leq 0.36577. \]
\end{restatable}

\begin{restatable}{theorem}{QnQlCk} \label{thm:QnQlC2k}
    For integers $k, \ell$ and $n$ such that $\ell \geq \log_2(2k)$, ~  $\de(Q_n, Q_\ell, C_{2k}) = 0.$
    For fixed integers $k,\ell$ and $n$ with $k \geq 4, k \neq 5$ and $2 \leq \ell < \log_2(2k) \leq n$, as well as $m: = \lceil \log_2(2k) \rceil -1$,  there is a positive constant $c_k$ such that
    \[ \binom{m}{\ell}\binom{n}{\ell}^{-1} \leq \de(Q_n, Q_\ell, C_{2k}) \leq c_k n^{-\frac{1}{16}}.\]
\end{restatable}

\begin{restatable}{theorem}{QnClQk} \label{thm:QnC2lQk}
    For integers $k, \ell$ and $n$ with $k\geq 2$ and sufficiently large $n$, we have
  \[ \max\left\{\left(1 -  \frac{1}{k} \right)\frac{(\ell-1)!}{2z_{\ell, \ell}}(1 - o(1)), 1- \frac{\ell}{k}\right\} \leq \de(Q_n, C_{2 \ell}, Q_k) \leq 1 -  \alpha \frac{\log k}{k2^k}, \]
  for a positive constant $\alpha$.
\end{restatable}

\begin{restatable}{theorem}{QnCsixCfour} \label{thm:QnC6C4}
    For $n$ sufficiently large, we have
   $ 0.03125 \leq  \de(Q_n, C_6, C_4) \leq 0.1625. $
\end{restatable}

After seeing our manuscript,  Balogh and Lidick\'y observed that  the upper bound in Theorem 6 could be improved  to 0.06875 using standard flag algebras, \cite{BL}.

\begin{restatable}{theorem}{QnClCk} \label{thm:QnC2lC2k}
    For fixed integers $k, \ell$ and $n$ such that $4 \leq k \leq 2^n, k \neq 5, 2 \leq \ell \leq 2^n$ and $\ell \neq k$, we have
    \[ \binom{n}{\ell}^{-1} \frac{2^{\ell-\ceil{\log_2(2 \ell)}}}{z_{\ell, \ell}} (1 - o(1)) \leq \de(Q_n, C_{2 \ell}, C_{2k}) \leq c_k n^{-\frac{1}{16}} . \]
\end{restatable}

In \Cref{lem:generalUB} we will prove the bound $\de(Q_n, T, H) \leq {\ex(Q_n, H)}/{||Q_n||}$. Using this bound and known bounds on the extremal number we obtain all upper bounds except for \Cref{thm:QnQlQk,thm:QnC4C6,thm:QnC6C4}.

\section{Notation and known results about $\ex(Q_n, H)$} \label{extremal-hypercube}

We often represent the vertices of $Q_n$  as binary vectors of length $n$, such that two vectors are adjacent if and only if the Hamming distance between them is one.  We will write these vectors simply as strings of symbols, such as $0010$, and we identify the vectors with these strings. Any copy of $Q_k$ in $Q_n$ can be represented by a vector with $n$ entries, where some $k$ of the entries are $\ordstar$, and all other entries are either $0$ or $1$. Then assigning the value $0$ or $1$ to each $\ordstar$ entry in every possible way yields every vertex in a copy of $Q_k$. Call such a representation of $Q_k$ a {\it  star representation}. The positions of stars are called {\it star positions}.
For example, a $Q_2$ can be written as $a \ordstar b \ordstar c$ where $a, b$ and $c$ are binary strings. Since edges correspond to copies of $Q_1$, they are represented by a  vector with one star.
We also write $\ones(a)$ to denote the number of ones in a binary string $a$.
The $i^{\rm th}$ {\it vertex layer} of $Q_n$ is the set of vertices with exactly $i$ entries equal to $1$, and the $i^{\rm th}$ {\it edge layer} is the set of edges with exactly $i$ entries equal to $1$ in their star representation. We denote the position of the star corresponding to an edge $e$ by $\ordstar(e)$.
For a positive integer $n$, let $[n]=\{1, \ldots, n\}$. For a subgraph $H$ of a graph $G$, we denote by $G-H$ a subgraph of $G$ with edge set $E(G)\setminus E(H)$.
\vspace{\baselineskip}

Erd\H{o}s~\cite{erdos2} was the first to ask how many edges a $C_{2k}$-free subgraph of the cube can contain. He conjectured that $\ex(Q_n, C_4)=\frac{1}{2} ||Q_n||(1+o(1))$.
Subsequently, there has been extensive effort devoted to determining the extremal numbers of cycles in the hypercube.
The best lower bound so far, $\ex(Q_n, C_4) \geq \frac{1}{2}\left( 1+n^{-1/2} \right) ||Q_n||$, is due to Brass {et al.}~\cite{brass} (valid when $n$ is a power of 4). The best upper bound due to Baber \cite{baber} is $\ex(Q_n, C_4) \leq 0.60318||Q_n|| (1+o(1))$.
Chung~\cite{chung} showed that $\ex(Q_n, C_6) \geq \frac{1}{4} ||Q_n||$.
She also proved for $k \geq 2$ that
\[ \ex(Q_n, C_{4k}) \leq c_k  n^{-\frac{1}{4}}  ||Q_n||. \]
Subsequently, Conder~\cite{conder} proved that $\ex(Q_n, C_6) \geq \frac{1}{3} ||Q_n||$, and the best known lower bound is due to Baber~\cite{baber}:  $\ex(Q_n, C_6) \leq 0.36577 ||Q_n|| (1+o(1))$.
Balogh {et al.}~\cite{balogh} proved a nearly identical though slightly worse bound, also using flag algebras.
F\"uredi and \"Ozkahya~\cite{furedi} extended this result  by showing in particular that  $\ex(Q_n, C_{4k+2}) = O(n^{-q_k} ||Q_n|| )$, where $q_k = 1/(2k+1)$ for $k\in \{3,5, 7\}$, and $q_k = 1/16 - 1/(16(k-1))$ for any other $k\geq 3$.
\vspace{\baselineskip}

To summarize,  $\ex(Q_n, C_{2k}) = \Theta(||Q_n||)$ when $k=2$ or $k=3$,
$\ex(Q_n, C_{2k}) = o(||Q_n||)$ for all  $k \geq 4$  and $k \neq 5$, and it remains unknown whether $\ex(Q_n, C_{10}) =  o(||Q_n||)$.
More specifically,  for  $k \geq 4, k \neq 5$
\begin{equation}\label{extremal-cycles} \ex(Q_n, C_{2k}) \leq c_k 2^{n-1} n^{\frac{15}{16}}.
\end{equation}

\vspace{\baselineskip}

In the course of investigating the extremal number $\ex(Q_n, C_{4k+2})$, F\"uredi and \"Ozkahya~\cite{furedi} also proved that $N(G,C_{4a}) \leq ||G|| O(n^{2a-2}) + O(2^nn^{2a-\frac{1}{2} + \frac{1}{2b}})$ for any $C_{4k+2}$-free graph $G$, where $G$ is a subgraph of $Q_n$, $k \geq 3$ and $4a + 4b = 4k+4$.
\vspace{\baselineskip}

Next, we discuss some known results for $\ex(Q_n, H)$.
Let $H$ be a subgraph of $Q_n$ and  $c(n, H)$ be the minimum size of a set $S$ of edges of $Q_n$ such that every copy of $H$ in $Q_n$ contains at least one edge from $S$.
Let $c(H)= \lim_{n \rightarrow \infty}  c(n, H) / ||Q_n||$ and so
$c(H)= 1 - \lim_{n \rightarrow \infty} {\ex (Q_n, H)}/{||Q_n||}$.
Alon, Krech and Szab\'o~\cite{AKS} showed that
\begin{equation}\label{aks}
    \alpha \frac{\log d}{d 2^{d}} \leq c(Q_d) \leq \frac{4}{d^2+2d +\epsilon},
\end{equation}
where $\epsilon\in \{0, 1\}$, $\epsilon \equiv d \pmod 2$, and $\alpha$ is a positive constant.
Offner \cite{O} proved that for a tree $T$ on a fixed number of edges $c(T)=1$, and for a $Q_d$-tree  $T'$ of cardinality~$k$,  $c(T')= c(Q_d)$.
Here, a $Q_d$-tree of cardinality $k$ is a union of $k$ copies $G_1, \ldots, G_k$ of $Q_d$ such that for any $i\geq 2$ there is $j<i$ such that $G_i\cap G_j$ is isomorphic to $Q_{d-1}$ and
$(G_i -  G_j)\cap (\cup _{\ell=1}^{i-1} G_\ell ) = \emptyset$.  This definition mimics the notion of a tree-width in the hypercube setting. Among other results, Offner \cite{O} proved a counting lemma:

\begin{lemma} \label{count-Offner}
    Let $\epsilon >0$ and $d\in \mathbb{N}$ be fixed, and let $n\rightarrow \infty$. If $H$ is a subgraph of $Q_n$ and $||H || \geq (1- c(Q_d) + \epsilon) ||Q_n||$, there are $\Omega(n^d2^n)$ copies of $Q_d$ in $H$.
\end{lemma}

Conlon~\cite{conlon}  extended these  results by showing that $\ex(Q_n, T) = o(||Q_n||)$ for a wider range of subgraphs $T \subseteq Q_n$, including all cycles $C_{2k}$ with $k \geq 4$ except for $C_{10}$, which is still an open case.
A subgraph $H$ of $Q_{\ell}$ is said to have a  {\it $k$-partite representation} if every edge of $H$ has exactly $k$ non-zero bits (stars and ones) and there is a function $\sigma: [\ell]\rightarrow [k]$ such that for each $e\in E(H)$,  $e=a_1 \cdots a_\ell$,  the image $\{\sigma(i_1), \ldots, \sigma(i_k)\}$ of the set of non-zero bits $\{a_{i_1}, \ldots, a_{i_k}\}$ of $e$ is $[k]$, i.e., distinct non-zero bits have distinct images.
One can also give a hypergraph formulation of this definition.
Specifically, for $e\in E(H)$, let $E_e$ be the set of non-zero positions of an edge $e$, for example if  $e=100\ordstar01$ then $E_e = \{1, 4, 6\}$.
Let $\cH= \cH(H)$ be a hypergraph on $\ell$ vertices  with hyperedge set $\{ E_e: e\in E(H)\}$.
Then $H$ has $k$-partite representation if $\cH$ is a $k$-uniform and $k$-partite hypergraph.

\begin{theorem} [Conlon \cite{conlon}]
    Let $H$ be a fixed subgraph of a hypercube. If, for some $k$, $H$ admits a $k$-partite representation, then $\ex(Q_n, H)= o(||Q_n||)$.
\end{theorem}

In addition to the extremal problem, it is natural to consider Ramsey-type statements about hypercubes.  In particular we say a graph $H$ is \emph{Ramsey} if for any $k$, there is an $n_0$ such that for any $n \ge n_0$ every edge coloring of $Q_n$ with $k$ colors contains a monochromatic copy of $H$ in one of the colors.  If a graph $H$ is not Ramsey, it is easy to see that $\ex(Q_n,H) = \Theta(||Q_n||)$.  Alon et al.~\cite{ARSV} gave a complete characterization of all graphs $H$ that are Ramsey.  In particular they proved that all even cycles of length at least $10$ have this property.

\section{Basic properties}

\begin{lemma} \label{lem:counting}
    For any natural $n\geq 3$ and $k\in [n]$,  $N(Q_n, Q_k) = \binom{n}{k} 2^{n-k}$
~   ${ \rm   and    } ~N(Q_n, C_6) = N(Q_3, C_6) \cdot N(Q_n, Q_3) = 16 \binom{n}{3} 2^{n-3}.$
Moreover, for any integers $n$ and $\ell$ with $2 \leq \ell \leq 2^{n-1}$,
    \[N(Q_n, C_{2 \ell}) = \sum_{k=\ceil{\log_2(2 \ell)}}^{\min\{\ell, n\}} \binom{n}{k} 2^{n-k} z_{k,\ell}.\]
    In particular, for a fixed $\ell$ and large $n$, we
    have $$N(Q_n, C_{2 \ell}) =   \binom{n}{\ell} 2^{n-\ell}z_{\ell, \ell}(1+o(1)).$$
\end{lemma}

\begin{proof}
    The first statement follows from choosing $k$ star positions in $\binom{n}{k}$ ways and filling the remaining $n-k$ positions with zeros and ones.
    To count $C_6$'s in $Q_n$, observe that each copy of $C_6$ belongs to a copy of $Q_3$ and this copy is determined uniquely. Thus, it is sufficient to count all $Q_3$'s in $Q_n$ and then count the number of $C_6$'s in $Q_3$.
    The former is done by the first statement of the lemma. For the latter, a $C_6$ in $Q_3$ can be formed by taking  vertices $000$, $111$, any two vertices in the first and in the second layer (this gives $9$ copies), or taking vertex $000$, all vertices of the first layer and any two vertices of the second layer (this gives $3$ copies),  then copies of $C_6$ can be formed by taking a vertex $111$, all vertices in the second layer and any two vertices in the first layer (this gives another $3$ copies), and finally there is one copy of $C_6$ using all vertices of first and second layers.  So, all together there are $16$ copies of $C_6$ in $Q_3$.

    Now, consider the edge set of a copy $C$ of $C_{2\ell}$ in $Q_n$ and let $k$ be the number of different star (flip) positions of these edges. Then $C$ is a subgraph of a unique $Q_k$, defined by those $k$ positions, where $k \leq n$. Also note that $k \leq \ell$ because each position that is flipped needs to be flipped again to get to the starting vertex, and $k \geq \ceil{\log_2(2 \ell)}$ because $Q_k$ has $2^k$ vertices and we need $2 \ell$ different vertices for $C$. Thus, in order to count the number of $C_{2\ell}$'s in $Q_n$, for each integer $k$, $ \ceil{\log_2(2 \ell)}\leq k \leq \min\{\ell, n\}$, choose $k$ star positions in $\binom{n}{k}$ ways, and fix the values for other positions in $2^{n-k}$ ways. Finally, consider the $2\ell$ binary vectors on the chosen $k$ positions so that they form a copy of $C_{2\ell}$, there are  $z_{k, \ell}$ ways to do this.

    The last equation follows because all other terms have order of magnitude of at most $n^{\ell-1} 2^n = o(n^\ell 2^n)$.
\end{proof}

\begin{lemma} \label{lem:generalUB}
Let $A$ and $B$ be graphs with $A \subseteq B\subseteq Q_n$ and the property that any copy of $A$ in $Q_n$ is in the same number of copies of $B$.  Then,  $\de(Q_n,B,H)  \le \de(Q_n,A,H).$
In particular, for any graphs $T, H \subseteq Q_n$ and integers $\ell, m$ such that $n \geq \ell \geq m \geq 1$,
$$\ex(Q_n, Q_\ell, H) \leq \binom{n-m}{\ell-m}\binom{\ell}{m}^{-1} 2^{m-\ell} \ex(Q_n, Q_m, H)$$
\text{ and }
$$\de(Q_n, T, H) \leq {\ex(Q_n, H)}/{||Q_n||}.$$
\end{lemma}
\begin{proof}
Let $G \subseteq Q_n$ be $H$-free and let it contain $\ex(Q_n,B,H)$ copies of $B$. Let every copy of $A$ be in exactly $M$ copies of $B$ in $Q_n$. Consider the sets $$X = \{(\tilde{A},\tilde{B}): \tilde{A} \subseteq \tilde{B} \subseteq G, \tilde{A} \cong A, \tilde{B} \cong B\} \mbox{ and }
Y= \{(\tilde{A},\tilde{B}): \tilde{A} \subseteq \tilde{B} \subseteq Q_n, \tilde{A} \cong A, \tilde{B} \cong B\}.$$
We have 
$|Y| = N(Q_n,B) \cdot N(B,A) = N(Q_n, A)M$, so  $M= N(Q_n,B) \cdot N(B,A)/N(Q_n, A)$. On the other hand  $|X| = N(G,B) \cdot N(B,A) \leq N(G, A)M$.
Thus  $N(G, B)/N(Q_n, B)  \leq N(G, A)/ N(Q_n, A) $.
Therefore 
\[
\de(Q_n, B, H) = \frac{N(G,B)}{N(Q_n, B)} \le \frac{N(G,A)}{N(Q_n, A)} \le \de(Q_n,A,H).
\]
Recall that $N(Q_n, Q_m) = \binom{n}{m}2^{n-m}$ and $N(Q_n, Q_\ell) = \binom{n}{\ell} 2^{n-\ell}$.
The last two statements of the lemma now follow by using $B=Q_{\ell}$ and $A=Q_m$ or $B=T$ and $A=K_2$, respectively.
\end{proof}

\begin{corollary} \label{cor:oIfoE}
    Let $T, H  \subseteq Q_n$ be fixed subgraphs of $Q_n$ and $\ex(Q_n, H) = o(||Q_n||)$. Then $\ex(Q_n, T, H) = o(N(Q_n, T))$ and thus $\de(Q_n, T, H) = o(1)$.
\end{corollary}

\begin{proof}
    By \Cref{lem:generalUB} we have $\de(Q_n, T, H) \leq \ex(Q_n, H)/||Q_n|| = o(||Q_n||)/||Q_n|| = o(1). $
\end{proof}

\begin{lemma} \label{lem:non-increasing-QnQlH}
    For any graph $H \subseteq Q_n$ and integer $\ell < n$ we have
    \[\de(Q_n, Q_\ell, H) \leq \de(Q_{n-1}, Q_\ell, H).\]
\end{lemma}

\begin{proof}
    Let $G$ be an $H$-free subgraph of $Q_n$  containing the largest number of copies of $Q_\ell$, i.e., $N(G, Q_{\ell}) = \ex(Q_n, Q_\ell, H)$. Consider triples $(Q, i, x)$ where $Q \cong Q_\ell, Q \subseteq G$, $Q$ contains no star in position $i$ and  the value in position $i$ is $x \in \{0,1\}$. We count these triples in two different ways.
    If we choose $i$ and $x$ there are at most $\ex(Q_{n-1}, Q_\ell, H)$ valid copies $Q$, and there are $2n$ possibilities to choose such $i$ and $x$.
    On the other hand, if we fix a $Q$, we must choose a position containing no star of $Q$, and then $x$ is already determined by $Q$. Thus there are $n-\ell$ ways to choose such $i$ and $x$. It follows that
    \[ (n-\ell) \cdot \ex(Q_n, Q_\ell, H) = (n-\ell) \cdot N(G, Q_\ell) \leq 2n \cdot \ex(Q_{n-1}, Q_\ell, H). \]
    Dividing both sides by $N(Q_n, Q_\ell)$ and rearranging the result gives us
    \begin{align*}
        \frac{(n-\ell) \cdot \ex(Q_n, Q_\ell, H)}{2^{n-\ell} \binom{n}{\ell}} &\leq \frac{2n \cdot \ex(Q_{n-1}, Q_\ell, H)}{2^{n-\ell} \binom{n}{\ell}} \\
        \iff \frac{\ex(Q_n, Q_\ell, H)}{2^{n-\ell} \binom{n}{\ell}} &\leq \frac{\ex(Q_{n-1}, Q_\ell, H)}{2^{n-1-\ell} \binom{n-1}{\ell}} \\
        \iff \de(Q_n, Q_\ell, H) &\leq \de(Q_{n-1}, Q_\ell, H).
    \end{align*}
\end{proof}

Recall that $z_{\ell, \ell}$ is the number of copies of $C_{2\ell}$ in $Q_\ell$ using exactly $\ell$ distinct  star positions on its edges. To bound this number, let $Z(\ell)$ be the set of words with elements from $\{1, \ldots, \ell\}$, where each word contains each symbol twice, but for $1 \le k < \ell$, no interval of $2k$ positions contains each symbol an even number of times.

\begin{lemma}\label{z-ll}
For any $\ell\geq 4$,  $z_{\ell, \ell} = |Z(\ell)|2^{\ell}/ 4\ell.$ In particular,
$z_{\ell, \ell} \leq (2\ell)! /4\ell.$
\end{lemma}

\begin{proof}
Let $\cC(\ell)$ be the set of  cycles of length $2\ell$ in $Q_\ell$ using $\ell$ star positions. For each such cycle $C$ fix the first edge $e_1$ arbitrarily, order the edges as $e_1, \ldots, e_{2\ell}$, let $s_i = \star(e_i)$ be the star positions of the edges, and let $s(C)=(s_1, \ldots, s_{2\ell})$. We call $s$ the {\it star list} of $C$.
Each symbol must appear at least twice in $s(C)$ since each flip of the coordinate should appear again. Since there are exactly $\ell$ symbols,  each appears exactly twice. Note that if $s(C)$ contains an interval of $2k$ positions containing each symbol an even number of times that is not $s(C)$ itself, then
the edges in $C$ corresponding to the edges in this interval  are consecutive edges of $C$ that form a cycle of length less than $2\ell$. Thus $s(C)$ has no such interval. On the other hand, each word from $Z(\ell)$ gives a star list of a cycle from $\cC$.

So, the problem of finding $|\cC(\ell)|$ is equivalent to finding  $|Z(\ell)|$.
Since we could choose elements of the first vertex of $C$  in $2^{\ell}$ ways and then could order the edges in $2\cdot 2\ell$ ways by shifts and change of direction, we see that
$z_{\ell, \ell} = |Z(\ell)| 2^{\ell}/ 4\ell.$

Next, we shall give the bounds on $|Z(\ell)|$. We call a word with the set of elements $\{1, \ldots, \ell\}$  {\it good} if each element appears exactly twice. Each word from $Z(\ell)$ is good. The total number of good words is  $(2\ell)!/ 2^{\ell}$. Indeed, note the element $1$ could be placed in its two positions in $\binom{2\ell}{2}=\ell(2\ell-1)/2$  ways, the element $2$ could be placed in $(2\ell-2)(2\ell -3)/2$  ways, and so on. Thus in particular $|Z(\ell)|\leq (2\ell)!/ 2^\ell$. This gives the desired upper bound on $z_{\ell, \ell}$.
\end{proof}

\section{Proofs of the main theorems}

\QnQlQk*

\begin{proof}
    {\bf Lower bound:} For the first expression, for $i=0, \ldots, k-1$, let $G_i$ be the union of $q^{\rm th}$ edge layers of $Q_n$ for all $q\not\equiv i \pmod k$. In  particular $G_i$  contains no copy of $Q_\ell$  as we need edges in $k$ consecutive layers for this.  Any copy $Q$ of $Q_\ell$ in $Q_n$ is contained in $k-\ell$ $G_i$'s.
    If $x_i$ is the number of copies of $Q_\ell$ in $G_i$, then we have $N(Q_n, Q_\ell) (k-\ell)= \sum_{i=0}^{k-1} x_i $.
    Thus there is $i \in \{0, \ldots, k-1\}$ such that $G_i$ contains $x_i\geq N(Q_n, Q_\ell) (k-\ell)/k$ copies of $Q_\ell$.
    \vspace{\baselineskip}

For the second expression in the lower bound we will use a construction of Alon, Krech and Szab\'o~\cite{AKS}.   Fix $n$ and $k$.  For $0 \le i < \left\lfloor (k+1)/2 \right\rfloor$ and $0 \le j < \left\lceil (k+1)/2 \right\rceil $, let $G(i,j)$ be the graph obtained by deleting the edge $l \ordstar r$ from $Q_n$ if and only if
    \[\ones(l) \equiv i \bmod{\left\lfloor \frac{k+1}{2} \right\rfloor} \quad \text{ and } \quad \ones(r) \equiv j \bmod{\left\lceil \frac{k+1}{2} \right\rceil }.\]

\noindent
    {\bf  Claim~} For any $i,j$,  $0 \le i < \left\lfloor (k+1)/2 \right\rfloor $ and  $0 \le j < \left\lceil (k+1)/2 \right\rceil$, $G=G(i,j)$ has no copies of $Q_k$.\\

     Note that if there are $m-1$ star positions in a vector, we can fill them with all zeros, $m-2$ zeros and one $1$, etc., producing $m$ consecutive integers as  number of ones in this vector and realising all modulo classes  modulo $m$.    If we look at the star representation of a $Q_k$, at least one edge using the $\left\lfloor (k+1)/2\right\rfloor^{\text{st}}$ star is not in $G$ since there are  $\left\lfloor (k+1)/2 \right\rfloor -1$  stars to the left and  $ k -  \left\lfloor (k+1)/2 \right\rfloor \geq  \left\lceil (k+1)/2 \right\rceil -1$ stars to the right, and thus some assignment of $0$'s and $1$'s to these stars gives an edge that meets the criteria for deletion.  For example, if $k= 7$ and we consider the $Q_k$ $010\ordstar 100 \ordstar\ordstar 001\ordstar 1110 \ordstar 101\ordstar 101\ordstar$, then by assigning $000$ to the first three stars and $110$ to the last three, the number of ones on the left is $0 \bmod 3$ and on the right is $0 \bmod 3$. So the edge $010{0} 100 {0}{0} 001\ordstar 1110 {1} 101{1} 101{0}$ of the $Q_k$ is not in $G(0,0)$. By assigning $110$ to the first three stars and $000$ to the last three, the number of ones on the left is $2 \bmod 3$ and on the right is $1 \bmod 3$. So the edge $010{1} 100 {1}{0} 001\ordstar 1110 {0} 101{0} 101{0}$ of the $Q_k$ is not in $G(2,1)$. This proves the claim.
    \vspace{\baselineskip}

Let $\G_k$ be the set of all graphs $G=G(i,j)$, and note $|\G_k| = \left\lfloor (k+1)/2 \right\rfloor \left\lceil (k+1)/2 \right\rceil$.  We shall try to average and see for a fixed copy $Q$ of $Q_\ell$, to how many $G$'s  from $\G_k$ it belongs to. Pick any $t$, $1\leq t\leq \ell$ and consider the $t^{\text{th}}$ star position of $Q$.  Some edge of $Q$ with a star in this position is not in $G(i,j)$ for at most $t(\ell-t+1)$ choices of $i$ and $j$. Indeed, on the one hand there are $t-1$ stars to the left of the $t^{\text{th}}$ star, so no matter how these are filled with zeros and ones, there are at most $t$ possible numbers of ones one can achieve to the left of the $t^{\text{th}}$ star. On the other hand, there are  $\ell - t$ stars of $Q$ to the right of the $t^{\text{th}}$ star. Thus no matter how these are filled with zeros and ones, there are at most $\ell-t+1$ possible numbers of ones one can achieve to the right of the $t^{\text{th}}$ position.
Summing over $t$, we have that there are at least $(|\G_k|- \sum_{t=1}^{\ell} t(\ell -t+1))$ graphs from $\G_k$ that contain $Q$.\\

Let  $x_G$ be the number of copies of $Q$ in $G\in \G_k$. We have the sum of $x_G$'s over all graphs in $\G_k$  is at least $N(Q_n, Q_\ell)(|\G_k|- \sum_{t=1}^{\ell} t(\ell -t+1))$.
Since $|\G_k|=  \left\lfloor (k+1)/2 \right\rfloor \left\lceil (k+1)/2 \right\rceil \geq k(k+2)/4$, by the pigeonhole principle
we have that there is a graph $G$ in $\G_k$ that is $Q_k$-free and has the following number of copies of $Q$:
$$N(G, Q_\ell) \geq   N(Q_n, Q_\ell)\left(|\G_k|- \sum_{t=1}^{\ell} t(\ell -t+1)\right)/|\G_k| \geq N(Q_n, Q_\ell) \left( 1 - \frac{4\binom{\ell+2}{3}}{k(k+2)} \right).$$

    \noindent  {\bf Upper bound:} Assume that $n = k$,  let $G$ be a $Q_k$-free subgraph of $Q_n=Q_k$. In particular, $G$ is a proper subgraph of $Q_n$, i.e., missed  at least one edge. Since an edge is in $\binom{n-1}{\ell-1}$ copies of $Q_\ell$ in $Q_n$ and $n = k$, the number of copies of $Q_\ell$ in $G$ is at most $\binom{k}{\ell} 2^{k-\ell} - \binom{k-1}{\ell-1} = \binom{n}{\ell} 2^{n-\ell} \cdot \left(1- 2^{\ell-k} \ell/k \right)$. For $n > k$, \Cref{lem:non-increasing-QnQlH} gives the extremal bound. Again, dividing by $N(Q_n, Q_\ell)$ concludes the proof for the first expression in the upper bound.

    For the second expression in the upper bound, we use Lemma \ref{lem:generalUB} and an upper bound
    $\ex(Q_n, Q_k) \leq (1- \alpha \log k/(k2^k))||Q_n||$,  that follows from (\ref{aks}).
    \vspace{\baselineskip}
\end{proof}

\QnCfourCsix*

\begin{proof}
    \noindent {\bf Lower bound: } We shall pick $C_4$'s, i.e., $Q_2$'s according to a parity condition described below. Then we define $G$ to be the union of these $C_4$'s and argue that $G$ is $C_6$-free. Specifically, pick a copy of $Q_2$ if its star representation 
    is $l \ordstar  \ordstar r$, the first star is in the odd position, the second star follows the first immediately,  $\ones(l) \equiv 0 \bmod 2$ and $\ones(r) \equiv 0 \bmod 2$. We call the vector $l$ the {\it prefix} and the vector $r$ the {\it suffix} of the selected $Q_2$.
%
%
    Observe that if a selected copy of $Q_2$ has stars in positions $2k+1$ and $2k+2$, then all edges of this $Q_2$ have stars in position $2k+1$ or in position $2k+2$.  Thus if two selected $Q_2$'s share an edge, say with a star in position $2k+2$, then they both belong to a same selected $Q_2$ with stars in positions $2k+1$. It follows that the selected $Q_2$'s are edge-disjoint.

    Let $G$ be a graph formed by the union of selected $Q_2$'s. We see that the number of $Q_2$'s in $G$ is at least the number of selected $Q_2$'s, that is at least (summing over the position $p$ of the first star and considering the parity of the prefixes and suffixes)

    $$\sum_{
    \begin{array} {l}
     p \in \{3, \ldots, n-2\}\\
     p \mbox{ is odd}
     \end{array}
      }   2^{p-1-1} 2^{n-p-1-1} + \underbrace{2^{n-2-1}}_{p = 1} \geq \frac{n}{2} 2^{n-4}.$$

    Next, we shall verify that $G$  does not contain any copies of $C_6$.  Assume otherwise, that there is a copy $C$ of $C_6$ in $G$.
    We see that vertices of $C$ have the same entries in some $n-3$ positions, the other three positions are $i_0, i_1,$ and $i_2$,  ~$i_0<i_1< i_2$.  
     Here are four examples of how those positions could be filled by vertices of $C$: 
    \begin{equation}\nonumber
        \begin{matrix}
            011\\
            111\\
            110\\
            010\\
            000\\
            001
        \end{matrix}
        ~~~~~~~~~~~~~~
        \begin{matrix}
            011\\
            111\\
            101\\
            001\\
            000\\
            010
        \end{matrix}
        ~~~~~~~~~~~~~~
        \begin{matrix}
            010\\
            110\\
            100\\
            101\\
            001\\
            011
        \end{matrix}
        ~~~~~~~~~~~~~~
        \begin{matrix}
            010\\
            110\\
            100\\
            101\\
            001\\
            000
        \end{matrix}
    \end{equation}

    For such a cycle $C$ we label the vertices $v_0, \dots, v_{5}$ so that  $v_5v_0$ and $v_iv_{i+1}$ are edges of $C$ for $0 \le i \le 4$.   Because all values in the vectors representing the vertices of $C$ are fixed except for the values in position $i_0, i_1$ and $i_2$, we will
    denote by  $\widetilde{v_j}$ a vector of length three with elements corresponding to the values of $v_j$ at positions $i_0, i_1, i_2$, in order.
    Let    $\widetilde{v_0}  = (\alpha, \beta, \gamma)$.  Note that the value in every  position $i_0, i_1, i_2$  is changed (flipped) exactly twice as we consider $v_0, v_1, \ldots, v_5, v_0$,   as otherwise we would use a vertex twice. So assume without loss of generality that $\widetilde{v_1} = (\alpha, 1-\beta, \gamma)$, i.e.\ that $i_1$ is flipped first, as otherwise we can just rotate the cycle such that the first edge has its star in position $i_1$.

    We first claim that neither $i_0$ and $i_1$, nor $i_1$ and $i_2$ are in positions $2k+1$ and $2k+2$ respectively for some $k$.  Assume otherwise, say $i_0 = 2k+1$ and $i_1 = 2k+2$. Once $\gamma$ has become $1-\gamma$, neither value in $i_0$ nor in $i_1$ can ever change again because the  parity of the suffix for the chosen $Q_2$'s would be different.  Thus our claim holds.

    Assume now that $\widetilde{v_2}= (\alpha, 1-\beta, 1-\gamma)$. Then $\widetilde{v_3} \neq (\alpha, \beta, 1-\gamma)$ because the parity of the suffix  of the $Q_2$ containing the edge $v_2v_3$ is no longer even. Thus $\widetilde{v_3} = (1-\alpha, 1-\beta, 1-\gamma)$. But now the value in $i_1$ cannot change because the parity of the suffix is still not even. The value in $i_2$ cannot change because now the parity of the prefix of the $Q_2$ containing $v_3v_4$ would no longer be even.

    The case with $\widetilde{v_2} = (1- \alpha, 1-\beta, \gamma)$ works symmetrically.  Thus such a $C$ cannot exist, and $G$ is $C_6$-free.

    By \Cref{lem:counting}, noting that $C_4=Q_2$,  we have $N(Q_n, Q_2) = n(n-1) 2^{n-3}$, so we get
    $\de(Q_n, C_4, C_6) \geq N(G, C_4)/N(Q_n, C_4) \geq n 2^{n-5}/(n(n-1) 2^{n-3}) \geq 0.25/n.$
    \vspace{\baselineskip}

    \noindent {\bf Upper bound:}
    For the upper bound, consider a $C_6$-free subgraph $G$ of $Q_n$. Note that any two copies of $C_4$ in $G$ are edge-disjoint since otherwise their union would contain a $C_6$. The total number of edges in $G$ is at most $\ex(Q_n, C_6) \leq 0.36577 n2^{n-1}$ by a result of Baber~\cite{baber}.
    Thus $N(G,C_4) \leq \frac{||G||}{4} \leq 0.36577 n2^{n-3}$.
    Dividing this number by $N(Q_n, C_4)$ we get, for large enough $n$,
    \[\de(Q_n, C_4, C_6) \leq \frac{N(G, C_4)}{N(Q_n, C_4)} \leq \frac{0.36577 n2^{n-3}}{n(n-1) 2^{n-3}} \leq 0.36578 \frac{1}{n}. \qedhere \]
\end{proof}

\QnClCsix*

\begin{proof}
    \noindent  {\bf Lower bound:} We use the $3$-coloring of Conder \cite{conder} to create a $C_6$-free subgraph $G$ of $Q_n$. For this, consider $G$ whose edges $e = l \ordstar r$ satisfy  $\ones(l) - \ones(r) \equiv 0 \pmod 3$. Then $G$ is $C_6$-free, see \cite{conder} for a proof of this.
    We now choose copies $Q$ of $Q_\ell$ in $Q_n$ by a condition depending on $\ell$, and show that each of them contains a copy $C(Q)$ of $C_{2\ell}$ which is a subgraph of $G$.
    For $\ell \geq 6$ we pick a $Q_\ell$ if and only if its star representation $p_0 \ordstar p_1 \ordstar \cdots \ordstar p_\ell$ satisfies $\ones(p_i) \equiv 0 \pmod 3$ for all $i \in \{0, \dots, \ell\}$.
    For example, if $\ell = 6$, we would pick $\ordstar 1101 \ordstar \ordstar 0 \ordstar \ordstar \ordstar 0$, but we would not pick $1 \ordstar \ordstar \ordstar 0 \ordstar \ordstar \ordstar 0$ since $\ones(p_0) \equiv 1 \pmod 3$.  Then the number of those $Q_\ell$'s is at least
    $\binom{n}{\ell} 2^{n-3\ell-2},$
    because we can choose $\ell$ stars out of the $n$ positions, fill all other positions but up to two left of each star and two to the right of the last star ($2\ell + 2$ positions) with $0$'s and $1$'s and use the reserved positions to force each $p_i$ to satisfy $\ones(p_i) \equiv 0 \pmod 3$.
    \vspace{\baselineskip}

   Let $\ell \geq 6$.  For a copy $Q$ of $Q_\ell$, we define $C(Q)$ as follows by giving the specific values in the star positions of $Q$:
    \[
    \begin{matrix}
        11100\cdots 00000\\
        11110\cdots 00000\\
        01110\cdots 00000\\
        01111\cdots 00000\\
        00111\cdots 00000\\
        \vdots \\
        00000\cdots 01110\\
        00000\cdots 01111\\
        00000\cdots 00111\\
        01000\cdots 00111\\
        01000\cdots 00110\\
        01000\cdots 00010\\
        01100\cdots 00010\\
        11100\cdots 00010\\
    \end{matrix}
    \]
    Starting with the first vertex and considering the vertices of the cycle in order, we see that until we reach the vertex corresponding to  $00000\cdots 00111$, each edge has either exactly three  $1$'s to the left and no $1$'s to the right of its star position, or no $1$'s to the left and three $1$'s to the right of its star position. As all parities between the star positions of $Q$ are $0 \pmod 3$, all those edges satisfy $\ones(l)-\ones(r) \equiv 0 \pmod 3$ and are thus in $G$. Also note that $00000\cdots 00111$ is vertex number $1 + (\ell-3) \cdot 2 = 2 \ell -5$ in our cycle.
    The last $5$ edges fulfill either the same parities as above, or have exactly one 1 to the left and one 1 to the right, and are thus also in $G$.
    \vspace{\baselineskip}

    For $\ell = 4$ or $\ell = 5$ we pick a $Q_\ell$ for $G$ if its star representation $p_0 \ordstar p_1 \ordstar \cdots \ordstar p_\ell$ satisfies $\ones(p_0) \equiv \ones(p_\ell) \equiv 0 \pmod 3$ and $\ones(p_1) \equiv \dots \equiv \ones(p_{\ell-1}) \equiv 1 \pmod 3$.
    For example, if $\ell = 4$, we would pick $\ordstar 11011 \ordstar \ordstar 1 \ordstar 0$, but we would not pick $1 \ordstar \ordstar \ordstar 01 \ordstar 0$ since $\ones(p_0) \equiv 1 \pmod 3$.  As before, we picked at least $\binom{n}{\ell} 2^{n-3\ell-2} \in \Omega\left( \binom{n}{\ell} 2^n \right)$ $Q_\ell$'s.
    For each chosen copy $Q$ of $Q_\ell$, define $C(Q)$ by assigning specific values to star positions of $Q$ as follows in the cases $\ell=4$ and $\ell=5$, respectively:
    \begin{equation}\nonumber
        \begin{matrix}
            0000\\
            1000\\
            1100\\
            1110\\
            1111\\
            0111\\
            0011\\
            0001
        \end{matrix}
        ~~~~~~~~~~~~~~
        \begin{matrix}
            00100\\
            01100\\
            01101\\
            01001\\
            11001\\
            11011\\
            10011\\
            10010\\
            10110\\
            00110\\
        \end{matrix}
    \end{equation}

    Manual checking of those cycles (using the modulo $3$ conditions mentioned above) yields that all edges are in $G$, and all positions of the corresponding $Q_\ell$ are used. Just as one example we check that $e = \ordstar000$ is indeed in $G$.  For this $e$, we have $l = p_0$ and $r = p_1 0 p_2 0 p_3 0 p_4$, which both satisfy $\ones(l) \equiv 0 \pmod 3$ and $\ones(r) \equiv 0 \pmod 3$ as $\ones(p_1) \equiv \ones(p_2) \equiv \ones(p_3) \equiv 1 \pmod 3$ and $\ones(p_4) \equiv 0 \pmod 3$.

    Since all $C(Q)$ use all star positions of their corresponding $Q$, and the values in non-star positions of different copies of $Q_\ell$ differ in some position, we know that $C(Q) \neq C(Q')$ if $Q \neq Q'$. Thus, we have
    $N(G, C_{2\ell}) \geq \binom{n}{\ell} 2^{n-3\ell-2}.$

    By \Cref{lem:counting} we have $N(Q_n, C_{2 \ell}) = \binom{n}{\ell} 2^{n - \ell} z_{\ell, \ell} (1+o(1))$, so dividing the bound on $N(G, C_{2 \ell})$ we just obtained by this number yields the lower bound.
    \vspace{\baselineskip}

    \noindent {\bf Upper bound:}  \Cref{lem:generalUB} and a result of Baber \cite{baber} that $\ex(Q_n, C_6)\leq 0.36577||Q_n||$ yields the upper bound.
\end{proof}

\QnQlCk*

\begin{proof}

    \noindent First note that $Q_\ell$ contains all even cycles of length at most $2^\ell$, and thus for $\ell \geq \log_2(2k)$ we have $\ex(Q_n, Q_\ell, C_{2k}) = 0$. Thus from now on assume that $k \geq 4, k \neq 5$ and $\ell < \log_2(2k)$.
    \vspace{\baselineskip}

    \noindent {\bf Lower bound:} Let $G$ be the union of all $Q_m$'s with stars in the first $m$ positions. By filling all other positions with 0 and 1 we see that $G$ has $2^{n-m}$ different $Q_m$'s which are pairwise vertex disjoint. As a $Q_m$ can only contain cycles with length at most $2^m < 2^{\log_2(2k)} = 2k$, $G$ does not contain any cycles $C_{2k}$.
    On the other hand, by \Cref{lem:counting}, $G$ contains $\binom{m}{\ell}2^{m-\ell} 2^{n-m}$ $Q_\ell$'s.
    \vspace{\baselineskip}

    \noindent {\bf Upper bound:}
    The upper bound follows again from Lemma \ref{lem:generalUB} and the bound from (\ref{extremal-cycles}):
    $\ex(Q_n, C_{2k}) \leq c_k \cdot 2^{n-1} \cdot n^{15/16}$, for integer $k\geq 4$, $k\neq 5$.
\end{proof}

\QnClQk*

\begin{proof}
    \noindent  {\bf Lower bound:} \\
    \noindent
   Let $G_j$ be  a subgraph $G$ of $Q_n$ that is a union of $i^{\rm th}$  edge layers, for all $i \equiv j \bmod k$.
    Then clearly $Q_n-G_j$ has no copies of $Q_k$ for any $j = 0, \ldots, k-1$ and $Q_n=G_0\cup \cdots \cup G_{k-1}$.
    Let $y=y(\ell)$ be the number of copies of $C_{2\ell}$ containing a given edge of $Q_n$.
   By counting the number $X$ of pairs $(C, e)$, where $C$ is a copy of  $C_{2\ell}$ containing the edge $e$, we see that $X = N(Q_n, C_{2\ell} )2\ell = ||Q_n||y$.
   Thus $y= 2\ell N(Q_n, C_{2\ell})/ ||Q_n||$.
   Since for some $j \in \{0, \ldots, k-1\}$, $||G_j||\leq ||Q_n||/k$, and each copy of $C_{2\ell}$ uses an even number of edges in each layer,  we have that the total number
  copies of $C_{2\ell}$ containing at least one edge in $G_j$ is at most
  $$\frac{||Q_n||}{k}  \frac{2\ell N(Q_n, C_{2\ell})}{ ||Q_n||} \frac{1}{2} = \frac{\ell}{k} N(Q_n,  C_{2\ell}).$$
  Therefore, at least $(1- \ell/k) N(Q_n,  C_{2\ell})$ copies of $C_{2\ell}$ are  in $Q_n-G_j$.
  This bound is non-trivial if $\ell<k$. \\

  For the other lower bound,  we shall again consider graphs $Q_n-G_j$ and count the number of copies of $C_{2\ell}$  completely contained in the edge-layers of $Q_n$.  Consider the $i^{\text{th}}$ edge layer of $Q_n$ and let $Y_i$ be the set of  copies of $C_{2\ell}$'s in this layer.
 Note that for any set $L$ of $\ell$ positions, there is  $C_L\in Y_i$, such that $C_L$ has star positions set $L$  and  such that restricted to these positions the vertices of $C_L$ are represented as the sums of a zero vector of length $\ell$ and
a binary vector corresponing to a vertex of  $C_{2\ell}$ that is in the first layer of $Q_\ell$.  The values of  the positions not in $L$ are fixed with  exactly $i-1$ ones.
 Let $y_i=|Y_i|$.  Then considering all copies of $C_{2\ell}$ that is in the first layer of $Q_\ell$, we have
 $$y_i \geq  \binom{n}{\ell}  \frac{\ell!}{2\ell} \binom{n-\ell}{i-1}.$$
  Here, the first term corresponds to the number of ways to choose the star position set $L$, the second term corresponds to the number of cycles of length $2\ell$ in the first edge-layer of $Q_\ell$, and finally the third term is the  number of ways to assign $i-1$ ones in positions not in $L$.
  Note that we are not over counting since for any two distinct sets $L$ and $L'$ of the $\ell$ star positions,  the respective cycles are different -- one is constant on $L'\setminus L$ and other changes the value in these positions.
   Consider $y= \sum_{i=0}^{n-1} y_i$ and let $j$ be an index such that $\sum_{i\equiv j \bmod k} y_i  \leq y/k$.
Then the number of $C_{2\ell}$'s in $Q_n-G_j$ is
   \begin{eqnarray*}
   N(Q_n-G_j, C_{2\ell}) & \geq  & \left(1-\frac{1}{k}\right) y \\
  & \geq &    \left(1-\frac{1}{k}\right)  \sum_{i=1}^{n-\ell+1}  \binom{n}{\ell}  \frac{\ell!}{2\ell}  \binom{n-\ell}{i-1}\\
  & = &  \left(1-\frac{1}{k}\right)\frac{\ell!}{2\ell} \binom{n}{\ell} 2^{n-\ell}. \\
   \end{eqnarray*}
    Recalling that $N(Q_n, C_{2\ell}) = (1+o(1)) \binom{n}{\ell} 2^{n-\ell} z_{\ell, \ell}$, we get the lower bound.\\
    \noindent

    \noindent {\bf Upper bound:}
    Using the upper bound on $\ex(Q_n, Q_k)$ that follows from (\ref{aks}) as well as \Cref{lem:generalUB}, the result follows.
\end{proof}

\QnCsixCfour*

\begin{proof}
    \noindent {\bf Lower bound:}  Recall that $z_{3,3} = N(Q_3, C_6)= 16$ by \Cref{lem:counting}, and each $Q_3$ contains exactly one $C_6$ using only edges in one edge layer.  Further, every $C_6$ is in exactly one $Q_3$, so the copies of $C_6$ in $Q_n$ that use edges in only one layer are 1/16 of the total number of $C_6$'s.    As in the proof of Theorem~\ref{thm:QnC2lQk}, for $j =0,1$, let $G_j$ be the subgraph of $Q_n$ that is a union of $i^{\rm th}$  edge layers, for all $i \equiv j \bmod 2$. Neither $G_0$ or $G_1$ contains a copy of $C_4$. Each of the $C_6$'s in $Q_n$ is in $G_0$ or $G_1$.  Thus one of these graphs is $C_4$-free and contains at least 1/32 of the total number of $C_6$'s.

    \vspace{\baselineskip}

    \noindent {\bf Upper bound:}  We consider a $C_4$-free subgraph $G$ of $Q_n$ and count copies of $C_6$ in $G$.
    The largest number of edges in a $C_4$-free subgraph of $Q_3$ is $9$, and there are at most three copies of $C_6$ in a $C_4$-free subgraph of $Q_3$ (realized by a subgraph on $9$ edges with three missing edges forming a matching of edges in three different directions, i.e., stars in different coordinates).
    Note that each edge in $G$ can be shared between at most $\binom{n-1}{2}$ $Q_3$'s.
    Let $X := |\{Q \subseteq Q_n \mid Q \cong Q_3, ||Q \cap G|| = 9 \}|$. Then
    \[ ||G|| \geq \frac{\sum_{Q \cong Q_3, Q \subseteq Q_n} ||G \cap Q||}{\binom{n-1}{2}} \geq \frac{9X + 0 \cdot \left( \binom{n}{3} 2^{n-3} - X \right)}{\binom{n-1}{2}}.  \]
    Thus $X \leq {||G|| \binom{n-1}{2}}/{9}$. On the other hand, if $Q$ is a copy of $Q_3$ in $G$ and $||Q \cap G|| \leq 8$ then at most one $C_6$ is in $Q \cap G$. So the number of $C_6$'s in $G$ is

    \begin{eqnarray*}
      N(G, C_6) & \leq &   3X + 1 \cdot \left( \binom{n}{3} 2^{n-3} - X \right)\\
       &\leq &\frac{2}{9} ||G|| \binom{n-1}{2} + \binom{n}{3} 2^{n-3} \\
        &\leq &\frac{2}{9} \cdot \ex(Q_n, C_4) \binom{n-1}{2} + \binom{n}{3} 2^{n-3} \\
        &\leq &\frac{2}{9} 0.60318n2^{n-1} \binom{n-1}{2} + \binom{n}{3} 2^{n-3} \\
        &\leq &2.60848 \cdot \binom{n}{3} 2^{n-3},
    \end{eqnarray*}
    where the upper bound for the extremal number $\ex(Q_n, C_4)$ is due to Baber \cite{baber}.
    By \Cref{lem:counting} we have $N(Q_n, C_6) = 16 \binom{n}{3} 2^{n-3}$, so dividing $N(G, C_6)$ by this expression gives us the upper bound on $\de(Q_n, C_6, C_4)$.
    \end{proof}

\QnClCk*

\begin{proof}
    \noindent  {\bf Lower bound:}
    Let $m := \ceil{\log_2(2 \ell)}$. As in the proof of the lower bound in \Cref{thm:QnQlC2k}, take the union of all $Q_m$'s with stars in the first $m$ positions. We see again that this results in $2^{n-m}$ different $Q_m$'s which are pairwise vertex disjoint.
    Taking one $C_{2 \ell}$ in each of them results in no other cycles, and thus yields a $C_{2k}$-free graph with $2^{n- \ceil{\log_2(2 \ell)}}$ $C_{2 \ell}$'s.

    By \Cref{lem:counting} we have $N(Q_n, C_{2 \ell}) = \binom{n}{\ell} 2^{n - \ell} z_{\ell, \ell} (1+o(1))$, so dividing $N(G, C_{2 \ell})$ by this number yields the lower bound.
    \vspace{\baselineskip}

    \noindent {\bf Upper bound:}
    As in the proof of the upper bound in \Cref{thm:QnQlC2k},  Lemma \ref{lem:generalUB} and (\ref{extremal-cycles}) imply the result.
\end{proof}

\section{Acknowledgements}
The research of the fourth author was supported by the grants Nemzeti Kutat\'asi, Fejleszt\'esi \'es Innov\'aci\'os Hivatal, NKFI 135800 and Institute for Basic Science, IBS-R029-C1.  The research of the first author was partially supported by DFG grant FKZ AX 93/2-1. The research of the third author was supported by a DAAD Award: Research Stays for University Academics and Scientists (Program 57381327).

\section{Alternative proof of a lower bound in Theorem \ref{thm:QnQlQk}}

    Let $G$ be the graph obtained by deleting an edge $l \ordstar r$ from $Q_n$ if and only if
    \[\ones(l) \equiv 0 \bmod{\floor{\frac{k-1}{2}}} \text{ and }  \ones(r) \equiv 0 \bmod{\ceil{\frac{k-1}{2}}}.\]

    {\bf Claim 1~} $G$ contains no copies of $Q_k$.\\
        If we look at the star representation of a $Q_k$, at least one edge using the $(\floor{(k-1)/2}+1)^{\text{st}}$ star is not in $G$ since there are $\floor{(k-1)/2}$ stars to the left and $\ceil{(k-1)/2}$ stars to the right, and thus some assignment of 0's and 1's to these stars will give an edge that meets the criteria for deletion.  For example, if $k= 7$ and we consider the $Q_k$ $010\ordstar 100 \ordstar \ordstar \hspace{.75mm} 001\ordstar 1010 \ordstar 101\ordstar 101\ordstar$, then by assigning $100$ to the first three stars and $110$ to the last three, the number of ones on each side of the middle star is a multiple of 4, so the edge $010{1} 100 {0}{0} 001\ordstar 1010 {1} 101{1} 101{0}$ of the $Q_k$ is not in $G$.  This proves Claim 1.
    \vspace{\baselineskip}

    {\bf Claim 2~} For fixed integers $a$ and $r$ with $0\leq a <r$, $ \sum_{k\geq 0} \binom{m}{a+rk} = \frac{2^m}{r} +o(2^m).$\\
        It is known, see for example Gould \cite{G} or Benjamin et al.\ \cite{BCK}, that
        for $\omega$ equal to the $r^{\text{th}}$ primitive root of unity,
        $$\sum_{k\geq 0} \binom{m}{a+rk} = \frac{1}{r} \sum_{j=0}^{r-1} \omega^{-ja} (1+\omega^j)^m.$$
        Then, in particular, we see that $$\sum_{k\geq 0} \binom{m}{a+rk}= \frac{1}{r}2^m + q(m,r,a),$$
        where $q(m, r, a) = \sum_{j=1}^{r-1} \omega^{-ja} (1+\omega^j)^m$, and by the triangle inequality
        $$|q(m, r, a)|  \leq r \cdot \max\{ |1+ w^j|: j\in \{1, \ldots, r-1\}\}^m \leq r \cdot (2-\epsilon)^m,$$
        for a positive $\epsilon$ depending on $r$. This proves Claim 2.
    \vspace{\baselineskip}

    We call a tuple $\alpha = (\alpha_0, \ldots, \alpha_\ell)$ of $\ell+1$ integers {\it sparse} if $\alpha_i\geq \sqrt{n}$ and
    $\alpha_0+\cdots +\alpha_\ell = n-\ell$. Note that the number of such sparse tuples is at least  $\binom{n-(\ell+1) \ceil{\sqrt{n}}}{\ell}  = \binom{n}{\ell}(1-o(1))$.  Let $X_\alpha$ be the set of all copies $Q$ of $Q_\ell$ in $Q_n$, such that $Q = a_0\ordstar a_1\ordstar \cdots \ordstar a_\ell$ and the length of each $a_i$ is $\alpha_i$. Thus  the size of $X_\alpha$ is $2^{n-\ell}$ for each $\alpha$.
    Let $X=\{X_\alpha: ~ \alpha \mbox{ is sparse} \}$. We see that
    \begin{equation} \label{sizeX}
        N(Q_\ell, Q_n)= |X|(1+o(1)).
    \end{equation}
    \vspace{\baselineskip}

    {\bf Claim 3 ~}
    Let $\alpha = (\alpha_0, \ldots, \alpha_\ell)$ be sparse.
    If  $Q\in X_{\alpha}$ is not a subgraph of $G$, i.e., $Q$ contains  a deleted edge then some value of $i$ with $1 \le i \le \ell$, we must have
    \begin{eqnarray}\label{cond}
        \ones(a_0) + \cdots + \ones(a_{i-1}) & \equiv & x \bmod{\floor{\frac{k-1}{2}}} \mbox{ and }\\
        \ones(a_{i}) + \cdots + \ones(a_\ell) &\equiv  & y \bmod{\ceil{\frac{k-1}{2}}}, \nonumber
    \end{eqnarray}
    where $x \in \{-i+1, -i+2, \ldots, 0\}$ and $y \in \{-\ell+i, -\ell + i+1, \ldots, 0\}$.
    \vspace{\baselineskip}

    Indeed, assume that  the above condition  (\ref{cond}) does not hold for $i$, say with first subcondition failing for $x$.
    Consider an edge $e$  of $Q$ with the star position corresponding to the $i^{\text{th}}$ star position of $Q$. That is, $e$ is obtained by assigning some zeros or ones to all the star positions of $Q$ except for the $i^{\text{th}}$.
    Let $x''$ be the number of ones assigned to the first $i-1$ star positions,  so $0\leq x'' \leq i-1$.
    Let  $x' = \ones(a_0) + \cdots + \ones(a_{i-1})$,  so $x'  \not\in  \{-i+1, -i+2, \ldots, 0\}$ modulo $\floor{(k-1)/2}$.
    The number of ones to the left of the star position of $e$ is  $x'+x'' \not\equiv 0 \bmod {\floor{(k-1)/2}}$. This implies that no edge of $Q$ with the $i^{\text{th}}$ star has been deleted.  So, if (\ref{cond}) fails for all $i$'s,  none of the edges of $Q$ are deleted. This proves Claim 3.
    \vspace{\baselineskip}

    Now, we shall upper bound $q_\alpha$, the total number of  $Q\in X_\alpha$ satisfying (\ref{cond}).
    For $i \in [\ell]$, the number of binary vectors $b_{\alpha, i}$ ($c_{\alpha, i}$) of length $\beta_{\alpha, i} = \alpha_0 + \dots + \alpha_{i-1}$ ($\gamma_{\alpha, i} = \alpha_i + \dots + \alpha_\ell$) with number of ones congruent to a specific value modulo $r=\floor{(k-1)/2}$ ($r'= \ceil{(k-1)/2}$) is $2^{\beta_{\alpha, i}}/r (1+o(1))$  (or $2^{\gamma_{\alpha, i}}/r' (1+o(1))$). Note also that $\beta_{\alpha, i} + \gamma_{\alpha, i} = n - \ell$.
    Since there are $i$ values for $x$ and $\ell-i+1$ values for $y$ in condition (\ref{cond}), we have
$$      q_\alpha\leq \sum_{i=1}^{\ell}\frac{i}{r} \frac{\ell-i+1}{r'} 2^ {n-\ell} (1+o(1))  \leq 2^{n-\ell} \frac{4}{k^2-2k} \cdot\binom{\ell+2}{3}(1+o(1)).$$
%
    Summing up over all sparse $\alpha$, we have that the number of $Q_\ell$'s that are in $X$ and that contain a deleted edge is at most $\binom{n}{\ell}2^{n-\ell} 4 \binom{\ell+2}{3}/ (k^2-2k)(1+o(1))$ because those $Q_\ell$'s have to satisfy~(\ref{cond}) by Claim 3.
    Since the number of copies of $Q_\ell$ that are not in $X$ is at most $o(N(Q_\ell, Q_n))= o(\binom{n}{\ell} 2^{n-\ell})$ by (\ref{sizeX}),
    we have that the number of copies of $Q_\ell$ in $G$ is at least $\binom{n}{\ell} 2^{n-\ell} (1 - \binom{\ell+2}{3}/ (k^2-2k)(1-o(1)).$
    By \Cref{lem:counting} we have $N(Q_n, Q_\ell) = \binom{n}{\ell} 2^{n-\ell}$. Dividing  the above  bounds by this quantity concludes the proof of the lower bound.
    \vspace{\baselineskip}

\section{Improved lower bound for Theorem \ref{thm:QnC2lC6}.}
    As in the previous section, we call a tuple $\alpha = (\alpha_0, \dots, \alpha_\ell)$ sparse if $\alpha_i \geq \sqrt{n}$ and $\alpha_0 + \dots + \alpha_\ell = n-\ell$. We also define $X_\alpha$ to be the set of all copies $Q$ of $Q_\ell$ in $Q_n$, such that $Q = p_0\ordstar p_1\ordstar \cdots \ordstar p_\ell$ and the length of each $p_i$ is $\alpha_i$.
     For a fixed sparse $\alpha$ we select all $Q =  p_0\ordstar p_1\ordstar \cdots \ordstar p_\ell \in X_\alpha$ satisfying $\ones(p_i) \equiv 0 \pmod 3$.
     Using the distribution of binomial coefficients in modulo classes as in the previous section, we see that the number of such $Q$'s is
   \[\prod_{i=0}^{\ell} \sum_{k\geq 0} \binom{\alpha_i}{0 + 3k}
    = \prod_{i=0}^{\ell} \frac{2^{\alpha_i}}{3} (1 - o(1))
    = \frac{2^{n-\ell}}{3^{\ell+1}} (1-o(1)).\]

    For $X=\{X_\alpha: ~ \alpha \mbox{ is sparse} \}$ we again have $|X| = \binom{n}{\ell} (1+o(1))$, and so the number of $Q =  p_0\ordstar p_1\ordstar \cdots \ordstar p_\ell \in X$ satisfying $\ones(p_i) \equiv 0 \pmod 3$ is $\binom{n}{\ell} 2^{n-\ell} {3^{-\ell-1}} (1-o(1))$.
    This gives an improvement of the lower bound by a multiplicative term $(4/3)^{\ell+1}$.
\end{document}